%% file: AvramannalsUB.tex
\documentclass[11pt]{article}
\usepackage{amsmath,amssymb,amstext,dsfont,fancyvrb,float,fontenc,graphicx,subfigure, theorem}

\title{\Large\bf  On D\"umbgen's exponentially modified Laplace continued fraction for
 Mill's ratio}
\author{\sc F. Avram}
\date{}

\cleardoublepage 
\usepackage[strict]{changepage}
\def\abstractname{Abstract -}   % <-----------------
\def\abstract{\begin{adjustwidth}{1cm}{1cm} \par    \footnotesize \noindent {\bf \abstractname}
\def\endabstract{ \end{adjustwidth} \smallskip }}
{\theorembodyfont{\itshape}}
{\theorembodyfont{\itshape}}
{\theorembodyfont{\itshape}}
{\theorembodyfont{\itshape}}
{\theorembodyfont{\itshape}}
{\theorembodyfont{\rm}}
{\theorembodyfont{\rm}}
{\theorembodyfont{\rm}}
{\theorembodyfont{\rm }}

\usepackage{color,verbatim}              % Need the color package
\usepackage{graphicx}
\usepackage{amsmath}
\usepackage{amssymb}
\usepackage{psfrag}
\usepackage{fullpage}
\usepackage{hyperref}
\input{layout.tex}

\definecolor{MyDarkBlue}{rgb}{0,0.08,0.50}
\definecolor{BrickRed}{rgb}{0.65,0.08,0}

\newtheorem{thm}{Theorem}
\newtheorem{lem}[thm]{Lemma}
\newtheorem{rem}[thm]{Remark}
\def\beL{\begin{lem}} \def\eeL{\end{lem}}
\def\beR{\begin{rem}} \def\eeR{\end{rem}}
\def\beC{\begin{cor}} \def\eeC{\end{cor}}

\begin{document}
\maketitle
\vskip 1.5em
 \begin{abstract}
 The approximation of the Gaussian cumulative distribution  $\Phi(x)$ or of the related Mills ratio
\beq \label{e:R} R(x):=\frac{1- \Phi(x)}{\f(x)}:=h(x)^{-1}\eeq where $\f(x)$
is the standard Gaussian density, and $h(x)$ is  its hazard rate,
 have a long history starting with Gauss and Laplace and continuing
nowadays \cite{gupta1965system,lee1992laplace,bryc2002uniform,brycaddendum,
winitzki2003computing,kouba2006inequalities}.
Below, we  improve an important family of bounds provided recently by D\"umbgen \cite{duembgen2010bounding}.
 \end{abstract}

\begin{keywords}
continued fraction, Mill's ratio, hazard rate
\end{keywords}

\begin{MSC}
   60-08, 11A55, 62Q05.
\end{MSC}

%%%%%%%%%%%%%%%%%%%%%%%%%%%%%%%%%%%%%%
\sec{Introduction}

A convenient starting point for the study of the Gaussian  Mill's ratio \eqref{e:R}  is  the first order
ODE
\beq \label{GODE} R'(x)= x \; R(x) -1, \quad R(0)= \sqrt{\pi/2}.\eeq

 The equation \eqref{GODE}  allows building  a Taylor expansion around $0,$ and
a formal Laurent expansion in negative powers at $\I$, due to Laplace:
 \bea R(x) =\frac 1{x}\le(1- \frac 1{x^2} + \frac {3!!}{x^4}-\frac {5!!}{x^6}+...\ri).\eea
 The latter is  divergent (though asymptotic in the sense
of Poincar\'{e}); however,  this problem  may be remedied by
considering  continued fractions, whose domain of convergence typically is larger than that of the series.
The passage from series to a continued fraction with denominators $1$ \cite[pg. 21]{lorentzen1992continued} may be achieved by using recursively the formula
\bea &&1+ \sum_{i=1}^\I (-1)^i a_i x^i \approx \\&&\le(1+ a_1 x + (a_1^2-a_2) x^2 + (a_1^3-2 a_1 a_2+ a_3) x^3 + (a_1^4- 3 a_1^2 a_2 +a_2^2 + 2 a_1 a_3 -a_4) x^4 ...\ri)^{-1}\eea
yielding
\bea R(x) &=&\frac 1{x}\le(1- \frac 1{x^2} + \frac {3!!}{x^4}-\frac {5!!}{x^6}+...\ri)=\frac 1{x}\frac 1{1+ \frac 1{x^2} - \frac {2}{x^4}+\frac {10}{x^6}-\frac {74}{x^8}...} \\&=&\frac 1{x}\frac 1{1+ \frac 1{x^2}(1 - \frac {2}{x^2}+\frac {10}{x^4}-\frac {74}{x^6}...)}=
\frac 1{x}\frac 1{1+ \frac {x^{-2}}{1 + 2 x^{-2}(1-\frac {3}{x^2}+\frac {21}{x^4}...)}}
\eea
\beq \label{Lcf} &&=
\frac 1{x}\frac 1{1+ \frac {x^{-2}}{1 +   \frac {2 x^{-2}} {1+\frac {3}{x^2}(1-\frac {4}{x^2}...)}}}=\frac 1{x}\frac 1{1+ \frac {x^{-2}}{1 +   \frac {2 x^{-2}} {1+\frac {3x^{-2}}{1+\frac {4x^{-2}}{1+...}}}}}=
\frac{1}{1+\frac{v}{1 + \frac{2 v}{1+\frac{3v}{1
+...}}}}, v=\frac{1}{x^2}.\eeq

This equation is related to the famous Laplace's
continued fraction \eqref{eq:La}, which yields alternating upper and lower bounds for Mill's ratio.
Tighter alternating  bounds were derived recently
by \cite{duembgen2010bounding}, by judicious modifications
of the last denominators. We propose further modifications which improve numerically on D\"umbgen's,
and seem (but are  not yet proved) to provide alternating  bounds as well.

%Note the last fraction with positive coefficients
%is a "S(tieltjes)-fraction".

{\bf Contents}.  A brief review of continued fractions is given in Section \ref{s:cf}.  Lee's and D\"umbgen's approaches to the Gaussian Mill's ratio  are reviewed in Section \ref{s:LD}. The new family of bounds is introduced and illustrated numerically in Section \ref{s:D}. In Section \ref{s:Prym} we discuss briefly  the possibility of extending this  approach for providing continued fraction bounds for other Pearson densities, like the Gamma density, which is of  interest in queueing, for example in asymptotic studies of retrial queues in the Halfin-Whitt regime.

\sec{A brief review of continued fractions \label{s:cf}}

{\bf Definitions}. Recall that a  continued fraction
$$b_0 +\mK_{1}%^\I
\le(\frac{a_k}{b_k}\ri)=b_0+\frac{a_1}{b_1+\frac{a_2}{b_2 + ...}}=b_0+\frac{a_1}{b_1+\frac{a_2}{b_2 + \mK_{3}\le(\frac{a_k}{b_k}\ri)}}$$
where $\mK_{n}(\frac{a_k}{b_k}):=
\fr{a_n}{b_n+\mK_{n+1}(\frac{a_k}{b_k})}$
is defined, when convergent, as the limit of the {\bf convergents} $R_n(x)=\frac{A_n(x)}{B_n(x)}$ obtained  by replacing $\mK_{n}\le(\frac{a_k}{b_k}\ri)$
with $0$.

$A_n$ and
$B_n$ satisfy both
the forward Wallis-Euler  recursion $x_{n}=b_n x_{n-1} + a_n x_{n-2}, n\geq 2,$ with respective
initial conditions $A_0=b_0, A_{1}=a_1+ b_0 b_1,$ and $B_0=1, B_{1}=b_1$, and may also be written as "continuant" determinants:
$$A_n=\det\le(\begin{array}{ccccccc} b_0&-1&\\a_1&b_1&-1&\\&a_2&b_2&-1&\\& &&\vdots&\\& &&&a_{n-1}&b_{n-1}&-1\\& &&&&a_{n}&b_{n}\end{array}\ri)$$

$$B_n=\det\le(\begin{array}{ccccccc} b_1&-1&\\a_2&b_2&-1&\\&a_3&b_3&-1&\\& &&\vdots&\\& &&&a_{n-1}&b_{n-1}&-1\\& &&&&a_{n}&b_{n}\end{array}\ri)$$

{\bf Transformations}. For any sequence $p_k \neq 0,  p_{0}=1$, the two fractions
\beq \label{e:tr} b_0+ \mK_{1}%^\I
\le(\frac{a_k}{b_k}\ri), \quad b_0+ \mK_{1}%^\I
\le( \frac{p_{k-1} p_k a_k}{p_k b_k}\ri)\eeq
 are equivalent (have the same convergents). Thus, appropriate choices of $p_k$ will simplify either the numerators or denominators, as desired.% \fn[4]{see %http://en.wikipedia.org/wiki/Generalized\_continued\_fraction %for more information}.

 {\bf Laplace's continued fraction}. Applying the   transformation \eqref{e:tr} to
 \eqref{Lcf} with $p_k=x$  and putting $a_k= k + \d_0(k)$, one arrives to
  Laplace's continued fraction
\beq \la{eq:La} \mK_{1}%^\I
\le(\frac{a_k}{x}\ri)=\frac{1}{x+\frac{1}{x + \frac{2}{x+\frac{3}{x
+...}}}}, \quad x>0,\eeq
which converges to $R(x)$ on $(0,\infty).$
{Another continued fraction associated to the Taylor expansion around $0$ was provided by Shenton.}
\beR Note that due to the repetition  of the numerator $1$, it is more natural here to start indexing $R_n$ by $n=0$, so that    the terminating fraction with numerator $n$ is denoted by $R_n$. Thus,
$$R_0=\frac{1}{x}, R_1=\frac{1}{x+\frac{1}{x}}, R_2=\frac{1}{x+\frac{1}{x + \frac{2}{x}}},
...$$
\eeR

\beR Another derivation of  Laplace's continued fraction may  be obtained, following Euler, by differentiating \eqref{GODE}, which
yields
\beq \label{ODErec} R^{(n)}(x)&=& B_n(x) \; R(x) -A_n(x) \Leftrightarrow R(x)=\frac{A_n(x)}{B_n(x)}+ \frac{R^{(n)}(x)}{B_{n}(x)},\\ A_{n+1}(x)&=& x A_n(x) + n A_{n-1}(x), \quad (A_0(x), A_{1}(x))=(0,1), \no \\
B_{n+1}(x)&=& x B_n(x) + n B_{n-1}(x), \quad (B_0(x), B_{1}(x))=(1,x), \no\eeq
\eeR
 see \cite{kouba2006inequalities}.

{\bf Modified continued fractions}. The computation of continued fractions is often achieved by
the backward recurrence
$$R_{m,n}=b_m+\frac{a_m}{R_{m+1,n}}, m=n-1,n-2,...,0$$
where $R_{m,n}=\mK_{m}^n
\le(\frac{a_k}{b_k}\ri)=b_m+\frac{a_{m+1}}{b_{m+1}+\frac{a_{m+2}}{b_{m+2} + ...\frac{a_{n}}{b_{n}}}}$.

The classic starting point is $R_{n,n}=b_n$, but the result may often be improved by starting with modified last denominators $R_{n,n}=\b_n=b_n + \g$, i.e. by using
$$R_n(x)=b_0+\frac{a_1}{b_1+} \frac{a_2}{b_2 +} ...+ \frac{a_{n-1}}{b_{n-1} +}
+ \frac{a_{n}}{b_n +\g }=\frac{(b_n+\g)A_{n-1}+ a_n A_{n-2}}{(b_n+\g) B_{n-1}+ a_n B_{n-2}}=\frac{A_{n}+\g A_{n-1}}{ B_{n}+ \g B_{n-1}},$$ \cite[Ch. 5.5]{lorentzen1992continued}.
Note that we have switched here to the one line convention of writing  continued fractions (in which the subcontinued fractions following a  $+$ or $-$ are realigned on the first line), and that parametrizing the last modified denominator
 by $\b_n=b_n + \g$  (developping around the
 "usual" continued fraction\; coefficient $b_n$)  simplifies some expressions.

The idea  is to   replace $b_n$ by an "ansatz" $\b_n$ approximating more closely  the exact value $R_{n,n}$ \cite{winitzki2003computing}. We will call this unknown value the "correct ansatz".

{\bf The limit ansatz}. Assuming $n$ is big enough so that $R_{n,k}$ varies slowly in $n$, one such approximation is the limit ansatz obtained by solving
\beq \label{e:Ani} R_{n,\infty}=b_n+\frac{a_n}{R_{n,\infty}}.\eeq

{\bf Alternating bounds}. As noticed already by Brouncker and Euler \cite{dudley1987some,khrushchev2008orthogonal}, the positivity of the continued fraction numerators and denominators implies that the convergents
 yield upper and lower bounds
 \beq \label{b1} R_2 \leq R_4 \leq ...R_{2 n} ...\leq  R ...\leq R_{2 n+1} \leq ...\leq R_3 \leq R_1, \eeq
 valid on the domain of convergence of the continued fraction.

   In particular, the convergents of the Laplace  continued fraction  yield bounds valid on $(0,\I)$ (see also \cite[Prop. 7]{kouba2006inequalities}).

 General error estimates
\beq \label{b3} |R-\R_{n}| <  \frac{n!}{B_{n} B_{n+1}} \eeq
 are also available\fn[4]{
 The relations \eqref{b1}, \eqref{b3} are consequences of the Euler identities
 \beq \label{b2}
&&R_{n}-\R_{n+1} = \frac{\prod_{i=1}^{n+1} (- a_i)}{B_{n} B_{n+1}}\\&&R_{2 n+1}-\R_{2 n-1} = -\frac{b_{2 n+1} \prod_{i=1}^{2 n}  a_i}{B_{2 n-1} B_{2n+1}}\no\\&&R_{2 n}-\R_{2 n-2} = \frac{b_{2 n} \prod_{i=1}^{2 n-1}  a_i}{B_{2 n-2} B_{2n}}\no\eeq}.

{\bf Uniform bounds on $[0,\I)$}. The  Laplace and Shenton  continued fractions are quite
efficient in their "natural domains", and this allows constructing efficient approximations  based on both. However, if one entertains the somewhat  academic wish to use a single
approximation valid on $[0,\infty)$, one must improve the quality of approximation at $0$ if the continued fraction is based on the series at $\I$, and viceversa.

Following \cite{lee1992laplace,duembgen2010bounding} in their  tribute to Laplace, we will consider here the continued fraction  based on the series at $\I$.
Two strategies suggest themselves:
\BEN \im  use {\bf rational two-point Pad\'e
approximants} $[k_0,k_{\infty}]$ fitting $k_0$ derivatives at $0$ and
$k_\I$ derivatives at $\I$ (these seem to have  been introduced by
Murphy and McCabe \cite{mccabe1976continued}).%\fn[4]{ \cite{mccabe1976continued} studied the similar problem of
%approximating Dawson's integral $\te{ -x^2} \int_0^x \te{ u^2} du$).
%For $k_0=k_\I$, this resulted in the continued fraction:
%$$\bar{r}(x)= \sqrt{\pi/2} e^{x^2/2}
%-R(x)=\frac{x}{1-\frac{x^2}{3 + \frac{2 x^2}{5-\frac{3 x^2}{7
%+...}}}}=\mK_{1}^\I \le(\frac{a_k}{b_k}\ri)$$ %where $a_k= (-1)^k k x^2 +
%\d_0(k) x$, $b_k=2 k+1$.}

Reasonable uniform approximations  are already obtained
with $k_\I=2, k_0 =1,2,...$ \cite{bryc2002uniform}, the simplest one with $k_0=1$ being $R_2(x) =%\frac{1}{\sqrt{2 \pi}}
\frac{(\pi -2) \sqrt{2 \pi} + x
(4-\pi)}{2(\pi-2) + x \sqrt{2 \pi} + x^2 (4-\pi)}=\frac{1}{x+\frac 1{\b_1(x)}}$, where $\b_1(x)=\frac{(\pi -2) \sqrt{2 \pi} + x
(4-\pi)}{2(\pi -2)  + x
(3-\pi)\sqrt{2 \pi}}$. The fit at $0$ is due here
to $\b_1(0)=\frac{ \sqrt{2 \pi} }{2  }$.

\im use cleverly chosen modified continued fractions, which, besides fitting at $0,$ achieve  possibly also a good approximation of the  "correct ansatz".

Applying the limit ansatz \eqref{e:Ani} to Laplace's continued fraction amounts to replacing the denominator
$x$ below the numerator $n$ by the "terminating denominator"
\beq \label{e:An} \b_n(x)=\frac x2 + \sqrt{\le(\frac x2\ri)^2+ n}.\eeq

An even better starting point $\b_n(x)=\frac x2 + \sqrt{(\frac x2)^2+ \g_n}$, with $\g_n=\b_n^2(0)$ defined in \eqref{e:b0} has been proposed by \cite{duembgen2010bounding}, by exploiting both the functional form of the limit ansatz, and the correct
behavior at $0$.

The simplest choice is chosing linear modifications
\bea \b_n(x)=\l_n x + \b_n(0).\eea

 D\"umbgen's  best results, confirmed here, are finally obtained with {\bf exponential type modifications}. Since there is no clear reason for that, we might call
this an "inspiration ansatz".

\EEN

\sec{Lee's and D\"umbgen's modified Laplace continued fractions \la{s:LD}}
One possible approach, taken by  \cite{lee1992laplace}, is to consider  {\bf doubly modified convergents}
\bea \label{modf} R_n(x)&=&b_0+\frac{a_1}{b_1+} \frac{a_2}{b_2 +} ...+ \frac{a_{n-1}}{b_{n-1} +}
+ \frac{\a}{b_n +\g }\\&=&\frac{(b_n+\g) A_{n-1}+ \a A_{n-2}}{(b_n+\g) B_{n-1}+ \a B_{n-2}}=\frac{A_{n}+\g A_{n-1}+(\a-a_n) A_{n-2}}{ B_{n}+ \g B_{n-1}+ (\a-a_n) B_{n-2}},\eea
 with both the last numerator and denominator modified, and where $b_n=b_n(x), \g=\g(x)$ may depend on $x$.

Consider the sign of the approximation error,  supposing, more generally, that $R(x)$ is the Mill's ratio of a density $f(x)$ satisfying
\beq \label{e:p}  f'(x)=-q(x) \; f(x), \eeq
where $q(x)$ is rational. Then, $R(x)$
satisfies  the first order differential  equation
\beq \label{ODEr} R'(x)= q(x) \; R(x) -1, \eeq
generalizing \eqref{GODE}.
%\beq \label{ODEgen}c(x) R'(x)= b(x) \; R(x) -a(x), \eeq
%where $a(x), b(x), c(x)$ are polynomials.
Then, if $\lim_{u \to \infty}\f(u) R_n(u)=0,$
the  approximation error
\bea \D_{n}(x)={\int_x^\I \f(u) du} - { \f(x)} R_n(x)\eea may be expressed as
\bea && \D_{n}(x)=\int_x^\I \f(u) du + \int_x ^\I(\f(u) R_n(u))' du  =\int_x ^\I \f(u)(1 + R_n(u)'-q(u) R_n(u)) du.\eea

While the sign of the last integral is hard to analyze, it is easier to control the sign of the integrand
\beq \label{op} \d_n(u)=-\frac{\D_n'(u)}{\f(u)}=1 + R_n'(u)-q(u) R_n(u):=(G R_n)(u)\eeq
where we note that $G$ is precisely the operator
defining our function of interest \eqref{ODEr}.
Providing upper/lower bounds may thus be achieved by  ensuring that  $\d_{n}(u)$ is negative/positive for all $u$ in the domain of convergence.

\beR Let us note also an expression for the second derivative:
\beq \label{op2} \d^{(2)}_n(u)=-\frac{\D_n''(u)}{\f(u)}= R_n''(u)-2 q(u) R_n'(u)+(q^2(u)-q'(u)) R_n(u)-q(u):=(G^{(2)} R_n)(u).\eeq \eeR

We turn now to  D\"umbgen's impressive "creative denominator modifications", whose numerical results suggest that Lee's double modifications are not necessary. The basis is again an analysis of the  sign of the derivative of the approximation error
$\d_{n}(u)$ defined in \eqref{op}, this times in terms of the modification $\b_n(x)$ \cite[Lem. 1]{duembgen2010bounding}.

\beL  \label{l:Due}
Let $\b(u)=\b_n(u)$ denote differentiable terminating  modified denominators for Laplace's continued fraction
$$\frac{1}{x+\frac{1}{x + \frac{2}{x+\frac{3}{x
+...\frac{n-1}{x+\frac{n}{\b_n(x)}}}}}}$$ of the Gaussian Mills ratio. Then:
\beq \label{op0} && \d_{n}(u)=\frac{(-1)^{n-1}(n)!}
{B_{n-1}(u)^2} \T G_n \b(u), \; \; \T G_n \b(u)=u \b(u) + \b'(u)+ n -\b^2(u). \eeq
% && \d_{n}^{(2)}(u)=\frac{(-1)^{n-1}(n)!}
%{B_{n-1}(u)^2} \T G_n^{(2)} \b(u),  \\&& \T G_n^{(2)} \b(u)=P_0(u, \b(u),\b'(u),\b''(u)) + P_1(u, \b(u),\b'(u),\b''(u)) \b(u) + P_2(u) \b^2(u)+ P_3(u) \b^3(u) \no \\&&P_0(x)=-2 \left(2 x^2-\beta ''(x) x+\left(x^2+1\right)
%   \beta '(x)^2+2 \left(x^2+2\right) \beta
%   '(x)+4\right), \no \\&& P_1(x)=4 x
 %  \left(x^2+3\right)+2 x \left(x^2+5\right) \beta
 %  '(x)-\left(x^2+1\right) \beta ''(x), P_2(x)=x^4+2
 %  x^2-5, P_3(x)=-x \left(x^2+5\right) \no\eeq
 %  \bea -2 \left(2 x^2-\beta ''(x) x+\left(x^2+1\right)
 %  \beta '(x)^2+2 \left(x^2+2\right) \beta
 %  '(x)+4\right),\\-6 x^4-39 x^2+\left(x^2+3\right) \beta
  % ''(x) x-2 \left(x^4+9 x^2+6\right) \beta
 %  '(x)-15,\\-x \left(x^4+5 x^2-12\right),x^4+9
 %  x^2+6 \eea
\eeL

{\bf Proof:} The equation may be established by induction. The operator $\T G_n \b(u)$, which provides the sign of $\d_n(u),$  is also
given on   \cite[pg. 7]{duembgen2010bounding}.% who works in terms
%of $h_n(x)=1/R_n(x)=x +\fr 1{x+ \fr{2}{x+ ...\fr{n}{\b_n(x)}}}$.

The next step towards producing   uniform bounds  valid on $[0,\infty)$ is to find conditions on the modified
 denominators  $\b_n(x)$ which give rise to a zero of the error at $0$.

 \beL \la{l:co} The equations for ensuring $\D_n(0)=0, \D_n'(0)=0, \D_n''(0)=0$ are linear in $\b_n(0), \b_n'(0), r_n:=\b_n''(0)/\b_n(0),$ with solutions:
 \beq \label{e:b0} &&\D_n(0)=0 \Leftrightarrow \b_n(0)=\sqrt{2}\frac{\Ga(n/2+1)}{\Ga(n/2+1/2)}\\
 &&\d_n(0)=0 \Leftrightarrow \b_n'(0) ={\b_n^2(0)-n}\\&&\d_n^{(2)}(0)=0 \Leftrightarrow r_n ={2(\b_n^2(0)-n-\frac 12)}.
 \eeq
  The constants $\b_n'(0)$ and $ r_n$ are positive.
 \eeL

\beR These formulas will produce two-point
 Pad\'e approximants, when applied to rational modifications $\b(x)$.   \eeR

 {\bf Proof:} The first formula  is obtained in \cite[(13),(14)]{duembgen2010bounding}, by  imposing recursively the condition
 $\D_n(0)=0 \Eq R_n(0)=\sqrt{\fr{\pi}2}$ on the successive errors
 $$\D_0=1-\F(x) -\frac{\f(x)}{\b_0(x)}, \D_1=1-\F(x) -\frac{\f(x)}{x+\frac{1}{\b_1(x)}}, $$
 yielding $\b_0(0)=\sqrt{\fr2{\pi}}, \b_1(0)=\sqrt{\fr{\pi}2}, ...$.
 In general, we may note that $R_{n}(\b_n,x)=R_{n-1}(x + \fr{n}{\b_{n}},x),$ yielding $\b_k(0)=\frac k{\b_{k-1}(0)}, k=1,2,...$.

 The second  formula follows from  \eqref{op0}. In
 \cite[Thm 2]{duembgen2010bounding}, it is presented as a favorite choice among several possible  linear modifications $\b_n(x)=\l_n x + \b_n(0)$, and a proof that it yields alternating bounds is offered, but without mention of the two-point Pad\'e connection.

 For the  third formula, which does not appear in
\cite{duembgen2010bounding}, it is enough to consider the case $\b_n(x)=\b(0)+  x(\b^2(0)-n) + x^2 \fr{\b''(0)}2$. A tedious computation yields that
$$\d^{(2)}(0)=0 \Eq R''(0)=R(0) \Lra \fr{\b''(0)}{\b(0)}=2(\b^2(0)-n-\fr 12).$$
Intriguingly, the same expression
  appears  in a different context on the bottom of  \cite[pg. 9]{duembgen2010bounding}. This topic deserves further attention, and we are investigating currently  whether the second order two-point Pad\'e condition  leads to
 alternating bounds, as suggested by our numerical results.

 The positivity follows from \cite[Lem. 3]{duembgen2010bounding}.

 \beQ These  results suggest the interesting problem
 of obtaining minimal solutions to the Riccatti inequations $\T G_n \b(u) \geq (\leq) 0, \for u \geq 0, $, with constraints $\b(0) =\sqrt{2}\frac{\Ga(n/2+1)}{\Ga(n/2+1/2)},$ which would provide an optimal modification of Laplace's continued fraction. \eeQ

 Next,   \cite[Lem. 2]{duembgen2010bounding} offers a simplified method of establishing alternating bounds, by replacing the requirement of
 strictly negative/positive derivatives $\D_n'(x)$
  by the weaker requirement of strictly negative/positive and unimodal derivatives, which is easier to impose. This idea is not exploited in our
  paper.

 Finally, \cite{duembgen2010bounding} raises  the dilemma of choosing between several possible functional forms for $\b_n(x)$.
\be \im
 The approximations $\b_n(x)=x+ \b_{n}(0)$ are not far from Lee's bound  $\b_n(x)=x+ \sqrt{n+1}$, since it may be  shown that $\b_n(0) \in (\sqrt{n+1/2},\sqrt{n+1}).$ However, both Lee's and D\"umbgen's linear approximations fare not so well numerically.
\im  \cite[Thm.1]{duembgen2010bounding} considers square root modifications, in which $n$ in the ansatz \eqref{e:An} is replaced by the constants
$\b_n^2(0)$ of \eqref{e:b0}.
\im \cite[Thm.2]{duembgen2010bounding} considers more
 general linear modifications $\b_n(x)=\l_n x+ \b_{n}(0)$, where $\l_n=\b_n'(0)$ is choosen
  to make also the first derivative
  $\D_n'(0)$ equal to $0$. By Lemma \ref{l:Due}, this requires solving  $ \b_n'(0)+ n -\b_n^2(0)=0,$ yieding $\l_n =\b^2_{n}(0)-n$
 \cite[Sec 5]{duembgen2010bounding}.

 \im Finally, \cite[Thm.3]{duembgen2010bounding} shows that the rational bounds may be considerably
 improved by using  {\bf exponential-type
 modifications} of the last denominators. \ee

\sec{Improved D\"umbgen's  exponentially modified continued fractions \label{s:D}}
   We have  implemented one step further D\"umbgen's idea of considering  exponentially modified continued fractions, by looking for exponential + linear modifications:\beq \b_n(x)=c_n x + \b_n(0) e^{- \sqrt{r_n} x}=(\l_n + r_n \b_n(0) ) x + \b_n(0) e^{- \sqrt{r_n} x},\eeq where the new constants $r_n$ are chosen to make the  second derivative
  $ \D_n''(0)$ equal to $0,$   which requires, cf. Lemma \ref{l:co}, \beq r_n =\fr{\b_n''(0)}{\b_n(0)}={2(\b_n^2(0)-n-\frac 12)}.\eeq

  The figures below compare the exponential, our improved exponential (practically indistinguishable from $0$), and the linear and square root  modifications. As expected, the square root (who does not fit any derivatives at $0$) loses always near $0$, but catches up with the linear later. The exponential modifications are always better,
  especially the new one proposed here. The maximum errors of the first four terms are $.00021, .000048, .000030,.000016$.

 \begin{figure}[h!]
\begin{center}
\includegraphics{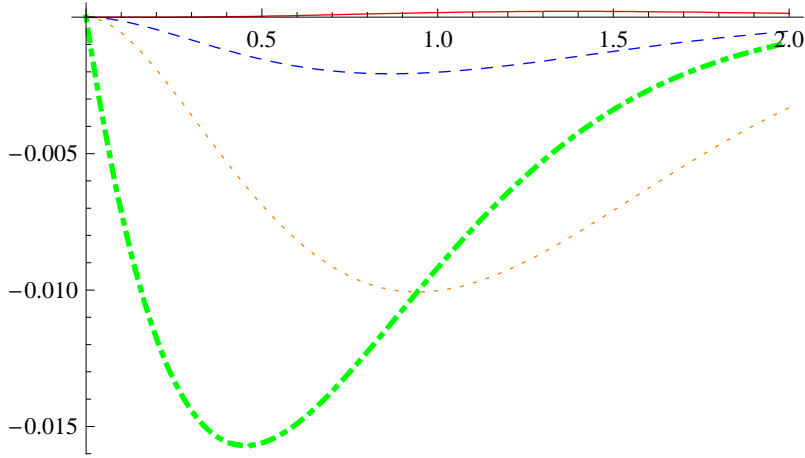}
\end{center}
\caption{Errors for D\"umbgen's bounds for $\D_0(x)$; blue, dashed: D\"umbgen's expo, red:second order expo, yellow, dotted: linear, green, dotdashed: square root}
\end{figure}

\begin{figure}[h!]
\begin{center}
\includegraphics{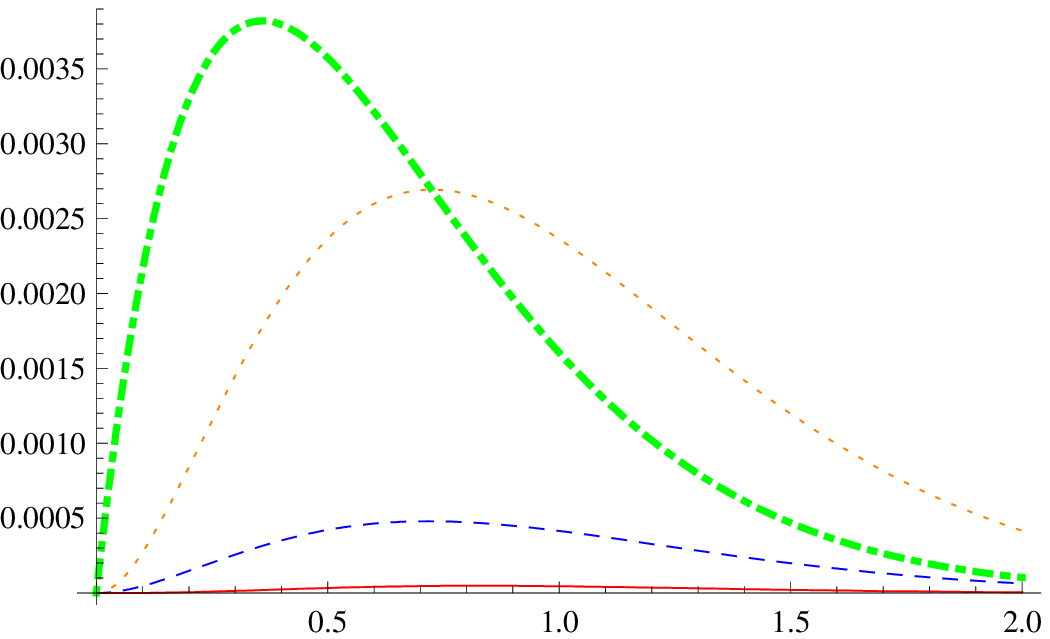}
\end{center}
\caption{Errors for D\"umbgen's bounds $\D_1(x)$; blue, dashed: D\"umbgen's expo, red:second order expo, yellow, dotted: linear, green, dotdashed: square root}
\end{figure}

\begin{figure}[h!]
\begin{center}
\includegraphics{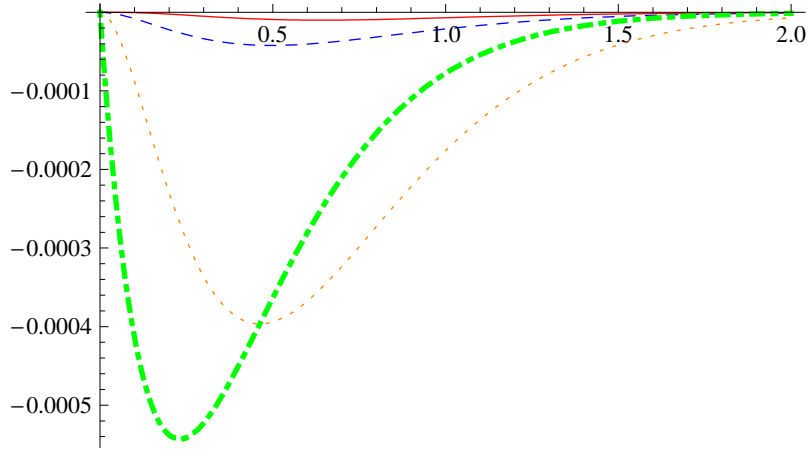}
\end{center}
\caption{Errors for D\"umbgen's bounds $\D_4(x)$; blue, dashed: D\"umbgen's expo; red:second order expo; yellow, dotted: linear; green, dotdashed: square root}
\end{figure}

\sec{Bounds for the
Gamma density Mills ratio/Prym's function \la{s:Prym}}

Besides  the normal,    bounds for  Mill's ratio of  other "Pearson distributions" (with connections to orthogonal
polynomials, etc...)  are also of great interest to probabilists.

The  Gamma density for example $\g(s,x)$   is  of special interest due to its appearance in many classic problems: the birthday paradox, Ramanujan's Q function,  Erlang loss probability, reliability, etc.
For the convenience of the reader interested in this problem, we summarize here some relevant information.

The Mills ratio $R(x)=R_s(x)$ for the
Gamma density $\g_s(x)$
 satisfies  the equation
\beq \label{ODEg} R'(x)= q(x) \; R(x) -1, \; q(x)=1+ \frac{1-s}{x}, \; R(\infty)=1 \quad (R(0)= 0, \text{ for } s <1)\eeq
and a  continued fraction for it was already developped in \cite{legendre1826traite}.
Note  the  integral representation:
 \beq \label{Rint} R_s(\l)=\l \int_0^\I(1+t)^s e^{- \l t} dt=\int_0^\I(1+\frac u\l)^s e^{- u} du.\eeq
   For integer $s,$ this may be easily derived by
 noting that the normalization of the Gamma density $\g_{k+1}(x)$ may be written as
 $\frac{s(s-1)...(s-k+1)}{\l^k}=\l \binom{s}{k}\int_0^\I t^k e^{-\l t} dt$ and summing for $k=0,1,...s$.
 For noninteger $s,$ see \cite{jagers1986continued}.

Some changes of  variables \cite[pg 143]{khovanskii1963application} put \eqref{ODEg} in the form of a homogeneous Riccati equation:
\beq \label{Ric} t^2 z'(t)- (1+ (1-s)t)  z(x) +z^2(t)=0,\eeq
from which the continuous fraction
\beq \label{L1} R(x)=\frac{x}{x+\frac{1-s}{1 + \frac{1}{x+\frac{2-s}{1
+...}}}}=\frac{x}{x}_{+\;}\frac{1-s}{1}_{+\;} \frac{1}{x}_{+\, \,}+\frac{2-s}{1}
{+...+}\frac{n}{x}_{+\, \,}\frac{n+1-s}{1}_{+...}\eeq
 may be obtained via a classic method of Lagrange \cite{de1776usage}. Cf \cite[(11.6)]{khovanskii1963application}, contracting the continuous fraction yields
 \beq \label{Laguerre}  e^x \int_x^\I u^{s-1} e^{-u} du=x^{s-1} R(x)= \frac{x^s}{x+1-s}_{+\;}\frac{s-1}{x+3-s}_{+\;} \frac{2(s-2)}{x+5-s}_
{+...+}\frac{n(s-n)}{x+ 2 n+1 -s}_{-...},\eeq
 a result which goes back to Laguerre
 \cite{laguerre1885reduction}. { A similar continued fraction\, expansion holds for the cumulative Gamma distribution:
 \beq \label{L2} x^{-s+1} e^x \int_0^x u^{s-1} e^{-u} du=\frac{x}{s-}\frac{s x}{1+s+x-} \frac{(1+s) x}{2+s+x-}...\frac{(n-1+s)x}{n+s+x-}...\eeq

 An equivalent  continued fraction used by  \cite{winitzki2003computing} is:
 \beq x^{1-s} e^x \int_x^\I u^{s-1} e^{-u} du=\frac{1}{1+}\frac{(1-s)v}{1+}\frac{v}{1+} \frac{(2-s)v}{1+}
\frac{2 v}{1+}{...}, \quad v=\frac 1 x.\eeq
This  generalizes  Laplace's  continued fraction. Indeed,  putting $s= \frac 12, u=\frac {v}{2}$ yields $\frac{1}{1+}\frac{u}{1+}\frac{2 u}{1+} \frac{(3 u}{1+}
{...}$, which is equivalent to Laplace's  Continued Fraction after substituting  $u=(\frac {1}{2 x})^2.$

The problem of providing bounds for the Gamma Mills ratio based on the continued fractions \eqref{L1}, \eqref{L2} has been considered by \cite{gupta1965system}. Several cases need to be distinguished, according to their difficulty:
 \BEN \im
 for $s \in(0,1]  $, the continued fraction approximations
continue to have positive coefficients, like in the Gaussian case (which corresponds in fact to $s=2,$ via a simple transformation). This case is thus straightforward \cite[(3.5)]{gupta1965system}.

\im for $s>1,$ fixed, the computation may be reduced to the case $s \in(0,1]$ by induction on the  integer part of $s$   \cite[(3.7)]{gupta1965system}.

\im the case $s \approx x \to \I$ is more subtle,   and it is precisely this case  that is of interest in queueing, for example in asymptotic studies of retrial queues in the Halfin-Whitt regime.
\EEN

 \section*{Acknowledgments} The author  thanks A. Baricz for having suggested further investigation of  D\"umbgen's approach.

 %%%%%%%%%%%%%%%%%%%%%%%%%%%%%%%%%%%%%%%%%
 %%%%  THE BIBLIOGRAPHY

 \bibliographystyle{alpha}
\bibliography{Mills}

 %%%%%%%%%%%%%%%%%%%%%%%%%%%%%%%%%%%%%%
 %%%%%%%%    THE AUTHORS

 { \footnotesize
\medskip
\medskip
 \vspace*{1mm}

\noindent {\it Florin Avram}\\
D\'epartement de Math\'ematiques, Universit\'e de Pau \\
Post address \\Avenue de l'Universit\'e - BP 1155, 64013 Pau Cedex, France.
    E-mail: {\tt florin.avram@univ-pau.fr}\\ \\

 \end{document}

%% file: layout.tex
    \newtheorem{Exa}{Example}
    \newtheorem{Exe}{Exercice}
\def\beXe{\begin{Exe}} \def\eeXe{\end{Exe}}

\def\eeD{\end{Def}} \def\beD{\begin{Def}}
\def\beXa{\begin{Exa}} \def\eeXa{\end{Exa}}

\newtheorem{Def}{Definition}
 \newtheorem{Qu}{Question}
\def\beQ{\begin{Qu}}
  \def\eeQ{\end{Qu}}
\def\beR{\begin{rem}} \def\eeR{\end{rem}}

\def\bep{\begin{pmatrix}}
  \def\eep{\end{pmatrix}}
  
  \def\Lra{\Longrightarrow}

    \def\BEN{\begin{enumerate}}  \def\BI{\begin{itemize}}
    \def\EEN{\end{enumerate}}   \def\EI{\end{itemize}} \def\im{\item}
       \def\sec{\section}
    
     \def\no{\noindent}

    \def\fr{\frac}

\newcommand{\beq}{\begin{eqnarray}} \newcommand{\eeq}{\end{eqnarray}}
\def\bea{\begin{eqnarray*}} \def\eea{\end{eqnarray*}}
\def\be{\begin{enumerate}} \def\ee{\end{enumerate}}
\def\bi{\begin{itemize}} \def\ei{\end{itemize}}
\def\im{\item} \def\I{\infty}
\def\D{\Delta} \def\a{\alpha } \def\d{\delta}
\def\b{\beta}  \def\g{\gamma}  

 \def\l{\lambda}

 \def\le{\left} \def\ri{\right}

\def\mbb{\mathbb}  
  
\def\R{\mathbb R}

\def\1{\mathop{\rm 1\!\!I}\nolimits}

 \def\Ga{\Gamma}

 \def\le{\left} \def\ri{\right} \def\I{\infty}
   \def\R{{\mathbb R}}

 \def\for{\forall} \def\T{\tilde}

 \def\fr{\frac}

\def\la{\label}

\def\Eq{\Leftrightarrow} \def\no{\nonumber}

\def\sec{\section}

\long\def\symbolfootnote[#1]#2{
\begingroup
\def\thefootnote{\fnsymbol{footnote}}\footnote[#1]{#2}
\endgroup}
\def\fn{\symbolfootnote}

\def\f{\phi}  \def\F{\Phi}

\def\Eq{\Leftrightarrow}

 \def\Ga{\Gamma}
\def\mK{{\mbb K}}
 
 \def\l{\lambda}

\def\ri{\right}

\newcommand{\ba}{\begin{array}}
\newcommand{\ea}{\end{array}}

{\begin{list}{}{\setlength{\rightmargin}{0cm}
                \setlength{\listparindent}{0cm}
                \settowidth{\labelwidth}{\mbox{#1}}
                \setlength{\leftmargin}{1.1\labelwidth}
                \setlength{\labelsep}{.1\labelwidth}}}%
{\end{list}}

\catcode`\@=11

\font\tenmsa=msam10 \font\sevenmsa=msam7 \font\fivemsa=msam5
\font\tenmsb=msbm10 \font\sevenmsb=msbm7 \font\fivemsb=msbm5
\newfam\msafam
\newfam\msbfam
\textfont\msafam=\tenmsa  \scriptfont\msafam=\sevenmsa
  \scriptscriptfont\msafam=\fivemsa
\textfont\msbfam=\tenmsb  \scriptfont\msbfam=\sevenmsb
  \scriptscriptfont\msbfam=\fivemsb

\def\Bbb{\ifmmode\let\next\Bbb@\else
 \def\next{\errmessage{Use \string\Bbb\space only in math mode}}\fi\next}
\def\Bbb@#1{{\Bbb@@{#1}}}
\def\Bbb@@#1{\fam\msbfam#1}